# On the Ermakov systems and nonlocal symmetries

Arunaye, I. F.

Department of Mathematics and Computer science
University of the West IndiesMona, Kingston 7. Jamaica.

E-mail: ifarunaye@yahoo.com; festus.arunaye@uwimona.edu.jm

## Abstract

Symmetry analysis of the Ermakov systems has attracted enormous treatments in recent times. In this paper we consider three classes of the Ermakov systems and obtain their nonlocal symmetries using a simple algebraic reduction process. We observed that these nonlocal symmetries are new to the literature.

## 1.0    Introduction

The Ermakov systems of second-order ordinary differential equations for which the invariant and the three generators of the algebra $sl(2,R)$ are well known in the literature (Ermakov [1], harin [2], Athorne [3], Goedert and Haas [4], Haas and Goedert [5], Simic [6], Leach and Karasu (Kalkanli) [7], Goodall and Leach [8]). The absence of sufficient numbers of Lie point symmetries for the Kepler problem in the context of complete symmetry group of dynamical systems brought to the fore the introduction of the nonlocal symmetries of dynamical systems by Krause [9]. Nucci [10] introduced her concept of reduction process combined with the Lie algorithm for obtaining the classical symmetries of differential equations to obtain the complete symmetry group of the Kepler problem according to Krause [9] and as well voided the earlier assertions of Krause [9] that these nonlocal symmetries could not be obtained by Lie algorithm. The Nucci [10] reduction process became so famous (even still) in the literature hitherto Arunaye [11] announced a simpler reduction process for reducing dynamical systems to systems of oscillator and equation of motion which admits Lie algorithm for the computation of their infinitesimal vector fields, meanwhile it is already established in Leach, Andriopoulos and Nucci [12] that the Ermanno-Bernoulli constants of dynamical systems are most suitable for reducing dynamical systems to systems of oscillator and equation of motion. Arunaye and

White [13] designated alternative constants for doing same; which were constructed from the Hamilton vector of dynamical systems. These alternative constants are equivalent to the Ermanno-Bernoulli constants in two dimensions and provide less cumbersome reduction variables in three-dimensions than the Ermanno-Bernoulli constants with less computational effort.

It is well known that the central feature of the Ermakov systems is their property of always having first integrals (Haas and Goedert [5], Simic [6]) and that this invariant plays a central role in the linearization of Ermakov systems Ray and (Reid [14], Haas and Goedert [5], Athorne [3]). The Kepler-Ermakov systems referred to the perturbations of the classical Kepler problem or an autonomous Ermakov system was investigated by Karasu (Kalkanli) and Yildirim [15] and they found that these systems are the usual Ermakov systems with frequency function depending on the dynamical variables. Leach, Karasu (Kalkanli), Nucci and Andriopoulos [7] supplemented the analysis of Karasu (Kalkanli) and Yildirim [15] by correcting some results; and also they carried out some investigations on the same dynamics in order to determine an equivalent transformation of the Kepler-Ermakov systems for the new time and rescaled radial distance so that in the discussion of the Kepler-Ermakov systems, suffices it to study its polar equivalent systems. Leach, Karasu (Kalkanli) [16] considered the Ermakov's superintegrable-toy for its nonlocal symmetries and asserted the insufficient Lie point symmetries and the unsuitable algebra $sl(2,R)$ for the complete specification of the system. It were further shown by the method of Nucci reduction process the five symmetries representation of the complete symmetry group; four of which are nonlocal symmetries and the algebra is the direct sum of a one-dimensional Abelian algebra and the semidirect sum of a solvable algebra with a two-dimensional Abelian algebra $[A_1 \oplus \{A_2 \oplus_s 2A_1\}]$ in the notation of the Mubarrakzyanov classification scheme Morozov [17], Mubarakzyanov [18, 19, 20]. This paper considered the reduction of three classes of the Ermakov systems by the method of Arunaye [11] for the nonlocal symmetries of same, the results are very interesting. The paper is organized as following. Section 2 dealt with the reduction of the three classes of Ermarkov systems to systems of oscillator and conservation law. Section 3 is devoted to the symmetry analysis and finally section 4 is concluding remarks.

## 1.1 Classes of Ermakov systems

The three classes of the Ermakov systems under consideration are given by

$$\ddot{x} + w^2(t)x = -\frac{x}{r^3}H + \frac{1}{x^3}f(\frac{y}{x}),$$
$$\ddot{y} + w^2(t)y = -\frac{y}{r^3}H + \frac{1}{y^3}g(\frac{y}{x}); \qquad (1.1)$$

$$\ddot{x}+w^2(t)x=\frac{1}{yx^2}f(\frac{y}{x}),$$

$$\ddot{y}+w^2(t)y=\frac{1}{xy^2}g(\frac{y}{x});\qquad(1.2)$$

and

$$\ddot{x}+w^2(t)x=\frac{1}{x^3},$$

$$\ddot{y}+w^2(t)y=\frac{1}{y^3};\qquad(1.3)$$

where $f$ and $g$ are arbitrary functions of their arguments, $H$ is a function of unspecified form of dependence upon $x$, $y$ and $r$ in which the Kepler part is to be found in the first term of the right sides of the system (1.1) (Leach and Karasu (Kalkanli) [16]). Systems (1.1), (1.2) and (1.3) are known as Kepler-Ermakov systems, generalized Ermakov systems and the Ermakov-toy systems respectively. In the symmetry analysis parlance, Leach and Karasu (Kalkanli) [16], Leach, Karasu (Kalkanli), Nucci and Andriopoulos [7] established the plane polar coordinates of (1.1.), (1.2) and (1.3) for their generalized symmetry analysis where $x=r\cos\theta$, $y=r\sin\theta$.

## 2.0   Reduction process of classes of Kepler-Ermakov systems

In the following the reduction process method for reducing dynamical systems to systems of oscillator and equation of motion reported in Arunaye [11] is applied to reduce these three classes of Ermakov systems (1.1), (1.2) and (1.3) to systems of oscillator and conservation law so that one is able to apply the Lie point symmetry analysis method to these classes of dynamical systems to obtain their generalized symmetries via their reduced systems.

## 2.1   Reduction process of the Kepler-Ermakov systems

The so called Kepler-Ermakov systems studied by Leach and Karasu (Kalkanli) [16] in polar system has radial and transversal components of the motion respectively given by

$$\ddot{r}-r\dot{\theta}^2=\frac{1}{r^3\cos\theta}h(\cot\theta)+\frac{1}{r^3}\{\sec^2\theta f(\tan\theta)+\cos ec^2\theta g(\tan\theta)\},\qquad(2.1)$$

$$r\ddot{\theta}+2\dot{r}\dot{\theta}=-\frac{1}{r^3}\{\sec^2\theta\tan\theta f(\tan\theta)-\cos ec^2\theta\cot\theta g(\tan\theta)\}.\qquad(2.2)$$

We note that the unspecified function $H$ in (1.1) and $h(\cot\theta)$ in (2.1) are related by $H = \frac{1}{4}Cr^3 - \frac{1}{r\cos\theta}h(\cos\theta)$ are C is an arbitrary constant (Leach and Karasu (Kalkanli) [16]).

Now from (2.2) we have

$$(r^4\dot\theta^2)^{\cdot} = 2\{\cos ec^2\theta \cot\theta g(\tan\theta) - \sec^2\theta \tan\theta f(\tan\theta)\}\dot\theta,$$

$$r^4\dot\theta^2 = L_\circ + 2\int\{\cos ec^2\theta \cot\theta g(\tan\theta) - \sec^2\theta \tan\theta f(\tan\theta)\}r^{-2}L dt.$$

i.e. $$L^2 = L_\circ + \alpha(\theta), \tag{2.3}$$

where $L_\circ$ is a constant, and $L = r^2\dot\theta$ defined the magnitude of the angular momentum of the motion which is not constant in this dynamics.

Now by setting $u = r^{-1}$; $\dot r = -Lu_\theta$; $\ddot r = -L^2 u^2 u_{\theta\theta}$ and substituting into (2.1), then we have

$$u_{\theta\theta} + [1 + \cos ec\,\theta h(\cot\theta) + \{\sec^2\theta f(\tan\theta) + \cos ec^2\theta g(\tan\theta)\}L^{-2}]u = 0.$$

i.e. $$u_{\theta\theta} + \omega^2(\theta)u = 0. \tag{2.4}$$

On taking $L_\circ = u_2$ and $u = u_1$ then (2.4) and (2.3) respectively become

$$u_{1,\theta\theta} + \omega^2(\theta)u_1 = 0, \tag{2.5}$$
$$u_{2,\theta} = 0,$$

where $\omega^2(\theta) = [1 + \cos ec\,\theta h(\cot\theta) + \{\sec^2\theta f(\tan\theta) + \cos ec^2\theta g(\tan\theta)\}L^{-2}]$. Thus (2.5) is the reduced system for the Kepler-Ermakov systems.

## 2.2 Reduction process of the Generalized Ermakov systems

The radial and transversal components of the motion for these dynamical systems are respectively given by

$$\ddot r - r\dot\theta^2 = \frac{1}{r^3}\{\sec^2\theta f(\tan\theta) + \cos ec^2\theta g(\tan\theta)\}, \tag{2.6}$$

$$r\ddot\theta + 2\dot r\dot\theta = -\frac{1}{r^3}\{\sec^2\theta \tan\theta f(\tan\theta) - \cos ec^2\theta \cot\theta g(\tan\theta)\}. \tag{2.7}$$

Following the same procedure as above, we obtain the similar reduced system for the dynamical systems (1.2) as

$$u_{1,\theta\theta} + \omega^2(\theta)u_1 = 0, \tag{2.8}$$
$$u_{2,\theta} = 0,$$

where $\omega^2(\theta) = [1 + \{\sec^2\theta f(\tan\theta) + \cos ec^2\theta g(\tan\theta)\}L^{-2}]$.

## 2.3 Reduction process of the Ermakov-Toy systems

The radial and transversal components of the motion for these dynamical systems are respectively known as

$$\ddot{r} - r\dot{\theta}^2 = \frac{1}{r^3}(\tan\theta + \cot\theta)^2, \tag{2.9}$$

$$r\ddot{\theta} + 2\dot{r}\dot{\theta} = -\frac{1}{2r^3}(\tan\theta - \cot\theta)', \tag{2.10}$$

where prime implies derivation with respect to $\theta$.

Similarly, the same procedure for reduction process above, reduces systems (1.2) to

$$u_{1,\theta\theta} + \omega^2(\theta)u_1 = 0, \tag{2.11}$$

$$u_{2,\theta} = 0,$$

where $\omega^2(\theta) = [1 + L^{-2}(\tan\theta + \cot\theta)^2]$.

## 3.1 Lie point symmetries of the reduced Ermkov-toy systems

We shall utilize (2.11) as hypothetical illustrative example of these three classes of Ermakov systems under investigation to present the nonlocal symmetries of the Ermakov-toy systems and note that it is easy to deduce the nonlocal symmetries for the Kepler-Ermakov systems and the generalized Ermakov systems.

Now we note that (2.11) is a system of $\theta$ dependent oscillator and a conservation law. The system (2.11) has nine Lie point symmetries (well known in the literature), they are

$$\begin{aligned}
\Gamma_1 &= \partial_\theta \ ; \ \Gamma_2 = \sigma^2 \sin 2\alpha \partial_\theta + u(\sigma\dot{\sigma}\sin 2\alpha + \cos 2\alpha)\partial_u \\
\Gamma_3 &= \sigma^2 \cos 2\alpha \partial_\theta + u(\sigma\dot{\sigma}\cos 2\alpha - \sin 2\alpha)\partial_u \ ; \ \Gamma_4 = \sigma \cos \alpha \partial_u ; \\
\Gamma_5 &= \sigma \sin \alpha \partial_u ; \ \Gamma_6 = \sigma^2 \partial_\theta + \sigma\dot{\sigma}u\partial_u ; \ \Gamma_7 = u\partial_u ; \\
\Gamma_8 &= \sigma u \sin \alpha \partial_\theta + u^2(\dot{\sigma}\sin \alpha + \sigma^{-1}\cos \alpha)\partial_u ; \\
\Gamma_9 &= \sigma u \cos \alpha \partial_\theta + u^2(\dot{\sigma}\sin \alpha - \sigma^{-1}\sin \alpha)\partial_u ;
\end{aligned} \tag{3.1}$$

where $\sigma$ satisfies the Pinney equation

$$\ddot{\sigma} + \omega^2(t)\sigma = \sigma^{-3}, \tag{3.2}$$

and $\alpha = \alpha(\theta)$ to be determined below. We note here that if $v$ and $\upsilon$ are two independent solutions of the oscillator in (2.11) for which $\sigma$ is as defined above then the solution of the Pinney equation satisfies

$$\sigma^2 = Av^2 + 2Bv\upsilon + C\upsilon^2, \ AC - B^2 = W^{-2} \tag{3.3}$$

where $A$, $B$ and $C$ are arbitrary real constants and $W := v\upsilon' - v'\upsilon$ is the wronskian where prime denotes derivation with respect to $\theta$. It is also well known that the Ermakov systems are famous as a result of the fact that they always posses the invariant

$$I_* = \tfrac{1}{2}[u^2/\sigma^2 + (\sigma\dot{u} - \dot{\sigma}u)^2]. \tag{3.4}$$

The transformation $u = r^{-1}$ implies $\partial_u = -r^{-2}\partial_r$ and $\partial_\theta = r^{-2}L\partial_t$, thus (3.4) becomes

$$I_* = \tfrac{1}{2}[\sigma^{-2}r^{-2} - (\sigma L^{-1}\dot{r} + \dot{\sigma}r^{-1})^2] \tag{3.5}$$

where $\sigma(t)$ is a function satisfying the Pinney equation. From the definition of the magnitude of the angular momentum $L$ we have

$$\dot{\theta} = r^{-2}L \Rightarrow \theta = \int r^{-2}L dt .$$

And one obtains

$$\alpha(\theta) = \int \frac{d\theta}{\sigma^2}. \tag{3.6}$$

We note that $\alpha(\theta)$ is a function of $t$ implicitly since $\theta$ is dependent on $t$.

### 3.2 Nonlocal symmetries of the Ermakov systems.

Substituting back for the original variables in the symmetries (3.1) of the Ermakov–toy systems we obtain the following generalized symmetries:

$$\begin{aligned}
&V_1 = r^{-2}L\partial_t \; ; \; V_2 = \sigma^2 r^{-2}L\sin 2\alpha\partial_t - r^{-3}(\sigma\dot{\sigma}\sin 2\alpha + \cos 2\alpha)\partial_r ; \\
&V_3 = \sigma^2 r^{-2}L\cos 2\alpha\partial_t - r^{-3}(\sigma\dot{\sigma}\cos 2\alpha - \sin 2\alpha)\partial_r ; \; V_4 = -\sigma r^{-2}\cos\alpha\partial_r ; \\
&V_5 = -\sigma r^{-2}\sin\alpha\partial_r ; \; V_6 = \sigma^2 r^{-2}L\partial_t - \sigma\dot{\sigma}r^{-3}\partial_r ; \; V_7 = -r^{-3}\partial_r ; \\
&V_8 = \sigma r^{-3}L\sin\alpha\partial_t - r^{-4}(\dot{\sigma}\sin\alpha + \sigma^{-1}\cos\alpha)\partial_r ; \\
&V_9 = \sigma r^{-3}L\cos\alpha\partial_t - r^{-4}(\dot{\sigma}\sin\alpha - \sigma^{-1}\sin\alpha)\partial_r , \; V_{10} = \partial_t ,
\end{aligned} \tag{3.7}$$

where $\sigma(t) = \left[\int \dfrac{d\eta}{\alpha(\eta)}\right]^{\frac{1}{2}}$ and $L^2 = L_\circ + \alpha(\theta)$ is obtained from (2.10) similarly as in (2.3). The $V_{10}$ is introduced as the symmetry responsible for the reduction of order by the change of independent variable from $t$ to $\theta$.

We observed that these generalized symmetries are absolutely different from those obtained in Leach and Karasu (Alkanli) [16], Karasu (Alkanli) and Yildirim [15] and Leach, Karasu (Alkanli), Nucci and Andriopoulos [7]. We also note that the point symmetries for time translation and special dilation respectively $V_1$ and $V_7$ have significant forms which are unparallel in the literature hitherto.

## 4.1  Conclusions

In the forgoing we are able to simply reduce the three classes of the Ermakov systems to systems of two equations - one time dependent oscillator another conservation law in a more simple and convenient method than the method of Nucci [10]. The application of Lie point symmetry algorithm is unique in that it produced the same nine point symmetries for each set of reduced systems however the backward transformation to obtain the nonlocal symmetries is much simpler too although the generalized symmetries obtained in this case may not be identical to those from Nucci reduction process consequence of the fact that nonlocal symmetries are infinite and have no unique algorithm for their general determination. We note the different manifestations of the nonlocal symmetries of these classes of Ermakov systems as with the realization of their complete symmetry groups, for which four are nonlocal symmetries and the algebra is the direct sum of a one-dimensional Abelian algebra and the semidirect sum of a solvable algebra with a two-dimensional Abelian algebra $[A_1 \oplus \{A_2 \oplus_s 2A_1\}]$ obtained in Leach, Karasu (Alkanli), Nucci and Andriopoulos [7]. The geometric implications as well as the uniqueness of the complete symmetry groups of these classes of Ermakov systems based on these new generalized (nonlocal) symmetries are subject for further discussions.

Reference:

[1] Ermakov V, Second-order differential equations: Condition of Complete integrability, Universita Izvestia Kiev ser III 30(1880), 1-25.

[2] Harin A.O., Second-order differential equations: Condition of Complete integrability, Appl. Anal. Discrete Math. 2 (2008), 123–145. *doi:10.2298/AADM0802123E*

[3] Athorne C, Kepler-Ermakov problems, J. Phys. A: Math. Gen. 24 (1991), L1385-L1389.

[4] Goedert J and Hass F, On the Lie symmetries of a class of generalized Ermakov systems, Phys. Lett. A 239 (1998), 348-352.

[5] Hass F and Goedert J, On the linearization of the generalized Ermakov systems. Phys. A: Math. Gen. 32 (1999), 2835-2844.

[6] Simic S S, A note on generalization of the Lewis invariant and the Ermakov systems, J. Phys. A: Math. Gen. 33 (2000), 5435-5447.

[7] Leach P G L , A Karasu (Kalkanli) , M C Nucci and K Andriopoulos Ermakov's Superintegrable Toy and Nonlocal Symmetries, SIGMA 1 (2005), Paper 018, 15 pages.

[8]  Goodall R and P G L Leach, Generalised Symmetries and the Ermakov-Lewis Invariant, J. Nonlinear Math. Phys. 12 (1) (2005), 15-26 (LETTER).

[9]  Krause J, On the complete symmetry group of the classical Kepler Systems, J. Math. Phys. 35 (11) (1994), 5734-5748.

[10] Nucci M C, The complete Kepler group can be derived by Lie group analysis.
    J. Math. Phys. 37 (4) (1996), 1772-1775.

[11] Arunaye F I On the reduction process of Nucci-REDUCE algorithm for computing nonlocal symmetries of dynamical systems: A case study of the Kepler problem, arXiv0709.1936 (September 2007), e-pub.

[12] Leach P G L, K Adriopoulos, and Nucci M C, The Ermanno-Bernoulli constants and representations of the complete symmetry group of the Kepler problem, J. Math. Phys. 44(9) (2003), 4090-4106.

[13] Arunaye F I and White H, On the Ermanno-Bernoulli and Quasi-Ermanno-Bernoulli constants for linearizing dynamical systems, Int. J. Appl. Math. and Informatics 1 (2) (2007), 55-60.

[14] Ray J R and Reid J L, More Exact Invariants for the time dependent harmonic Oscillator, Phys. Lett. A 71 (1979), 317-318.

[15] Karasu (Kalkanli) A and Yildirim H, On the Lie symmetries of Kepler-Ermakov systems, J. Nonlinear Math. Phys. 9 (4) (2002), 475-482.

[16] Leach P G L and Karasu (Kalkanli), The Lie algebra sl(2,R) and so-called Kepler-Ermakov systems, J. Nonlinear Math. Phys.11 (2) (2004), 269-275.

[17] Morozov V V, Classification of six-dimensional nilpotent Lie algebras, Izv. Vys. Uchebn. Zaved Matematika, 4(5) (1958), 161-171.

[18] Mubarakzyanov G M, On solvable Lie algebras, Izv. Vys. Uchebn. Zaved Matematika, 1(32) (1963a), 14-123.

[19] Mubarakzyanov G M, Classification of real structures of five-dimensional Lie algebras, Izv. Vys. Uchebn. Zaved Matematika, 1(32) (1963b),114-123.

[20] Mubarakzyanov G M, Classification of solvable Lie algebras of sixth-order with a non-nilpotent Lie algebras, Izv. Vys. Uchebn. Zaved Matematika, 3(34) (1963c), 99-106.